\documentclass[journal]{IEEEtran}

\usepackage{cite}
\usepackage{graphicx}
\usepackage{color}
\usepackage{amsmath,amssymb,amsfonts}
\usepackage{multirow}
\usepackage{url}
\usepackage[version=3]{mhchem}

\def\Vec#1{\mbox{\boldmath $#1$}}

\begin{document}
\bstctlcite{IEEEexample:BSTcontrol}
	
\title{Analysis of Heterogeneity of Pneumothorax-associated Deformation using Model-based Registration}

\author{Megumi~Nakao,~\IEEEmembership{Member,~IEEE,}
	Kotaro~Kobayashi,
	Junko~Tokuno,
	Toyofumi~Chen-Yoshikawa,
	Hiroshi~Date,
	and~Tetsuya~Matsuda,~\IEEEmembership{Member,~IEEE}% <-this % stops a space
	
	\thanks{M. Nakao, and T. Matsuda are with Graduate School of Informatics, Kyoto University, Kyoto, Yoshida-Honmachi, Sakyo, Kyoto, JAPAN. e-mail: megumi@i.kyoto-u.ac.jp}% <-this % stops a space
	\thanks{Dept. of Thoracic Surgery, Kyoto University Hospital, 54 Kawaharacho, Shogoin, Sakyo, Kyoto, 606-8507, Japan.}% <-this % stops a space
	\thanks{Dept. of Thoracic Surgery, Nagoya University Hospital,65 Tsurumai-cho, Showa-ku, Nagoya, 466-8560, Japan.}
	
	\thanks{Manuscript received October 20, 2020;}}

% The paper headers
\markboth{Journal of \LaTeX\ Class Files, October~2020}%
{M. Nakao \MakeLowercase{\textit{et al.}}: Analysis of Heterogeneity of Pneumothorax-associated Deformation}

% make the title area
\maketitle

\begin{abstract}
Recent advances in imaging techniques have enabled us to visualize lung tumors or nodules in early-stage cancer. However, the positions of nodules can change because of intraoperative lung deflation, and the modeling of pneumothorax-associated deformation remains a challenging issue for intraoperative tumor localization. In this study, we introduce spatial and geometric analysis methods for inflated/deflated lungs and discuss heterogeneity in pneumothorax-associated deformation. Contrast-enhanced CT images simulating intraoperative conditions were acquired from live Beagle dogs. Deformable mesh registration techniques were designed to map the surface and subsurface tissues of lung lobes. The developed framework addressed local mismatches of bronchial tree structures and achieved stable registration with a Hausdorff distance of less than 1 mm and a target registration error of less than 5 mm. Our results show that the strain of lung parenchyma was 35\% higher than that of bronchi, and that subsurface deformation in the deflated lung is heterogeneous. 
\end{abstract}


\IEEEpeerreviewmaketitle

\section{Introduction}

% Background
Recent advances in medical imaging devices have increased the chances of finding early-stage cancer, metastatic lung tumors, and benign nodules. Although the positions of lung nodules can be identified preoperatively on computed tomography (CT) images, the lungs deform or collapse during surgery, which prevents accurate tumor localization and increases the difficulty of resection procedures. In the last decade, image-based and model-based intraoperative guidance have been explored in laparoscopic surgery \cite{Koo17, Heiselman18}; however, pneumothorax-associated deformation has not been targeted because of a lack of clinical datasets and the complex properties of lung tissue. The analysis and modeling of pneumothorax-associated deformation and its application to intraoperative tumor localization therefore remains a challenging area of research for biomedical engineering.

% Pneumothorax deformation and technical issues
In the image-based analysis of lung modeling, respiratory motion has been extensively investigated \cite{Ehrhardt11, Naini12, Ilegbusi12, Fuerst15, Wilm16, Jud17}. However, because of difficulties in both image registration and regularization, few studies have modelled the pneumothorax-associated deformation that can occur between the preoperative and intraoperative lung states. When deformable registration is applied to CT images of pneumothorax, two technical issues should be addressed: the CT intensity shift and the major deformation with both rotational and topological changes. The CT intensity shift derives from the atelectasis state, in which there are variations in air content within the lungs; this changes pixel values and severely lowers image contrast, seriously impairing image-based registration. Lungs are very soft and their deformation is accompanied by significant volume reduction inside the thoracic cavity, which results in topological changes in the CT image space. The mechanism of pneumothorax-associated deformation is complex and not mathematically understood, except through simulation studies of animal lungs \cite{Naini12, Naini11}.  

% Related/previous studies
Deformable image registration techniques \cite{Sotiras13, Oh17, Faisal05} have been used for numerical analyses of organ and soft tissue deformation. Intraoperative deformation derived from changes in internal pressure, patient posture, and tool manipulation has been studied in respect to intraoperative guidance systems \cite{Kenngott14, Wagner16, Nickel18, Koo17, Ruhaak17, Nakao10, Nakao14, Gunay16, Tokuno20}, and studies have reported results of deformable registration designed for guidance of thoracoscopic surgery. Uneri et al. presented a deformable registration framework for pneumothorax cone beam CT (CBCT) data of animal lungs \cite{Uneri13}, and an integrated framework of model-based and image-based registration was used to address CT intensity shifts; however, the main focus of the work was registration of bronchial junctions. Alvarez et al. analyzed CBCT images of lungs deformed because of posture differences \cite{Pablo18}, and Nakao et al. designed a deformable mesh registration (DMR) framework for registering pneumothorax-associated lung shapes and statistically investigated lung surface deformation in animal CT datasets \cite{Nakao19}. More recently, the deflated lung surfaces of patients with pneumothorax were evaluated using hyperelastic finite-element models \cite{Jeanne20}. However, to the best of our knowledge, no study has reported registration results for whole lung structures or analyzed their spatial displacements, including their surfaces and subsurface structures. Therefore, we performed data acquisition and detailed numerical analysis with the aim of generating a statistical formulation of pneumothorax-associated deformation.

% Purpose
The purpose of this study was to perform a spatial and geometric analysis of pneumothorax-associated deformation to better understand the intraoperative state of the lung. Specifically, we sought to determine the heterogeneity of 3D displacements as evidenced by targeting of the lung’s surface and subsurface, including bronchial structures. In this analysis, we used pairs of contrast CT datasets acquired at two bronchial pressures (designated as inflated and deflated states) from 11 live Beagle dogs. To address the difficulties of mapping pneumothorax CT images with topological changes and CT intensity shifts, we designed a DMR framework for mixed data structures including the lobe surfaces and part of the centerline of the bronchi. We used the spatial maps obtained from the inflated/deflated lungs to investigate differences in strain (a measure of the displacement between particles in the body relative to a reference length) on the lung’s parenchyma and bronchial structures. 

% Contribution
The contributions of this paper are summarized as follows.
\begin{itemize}
    \item a complete model-based deformable registration solution for mapping inflated/deflated lungs
    \item an objective function design that addresses local mismatches of subsurface structures derived from CT intensity shifts 
    \item visualization of deformation fields including per-lobe contractions and rotations
    \item numerical analysis of heterogeneity in pneumothorax-associated deformation of whole lungs
\end{itemize} 
The registered models and findings can be directly used for statistical modeling of lung deflation. The potential clinical applications include localization of deep tumors or nodules and development of intraoperative visual guidance for thoracoscopic surgery.

\section{Material and Methods }
\subsection{Measurements and preprocessing}
To analyze spatial displacements in collapsed lungs, we used CT images acquired from 11 live Beagle dogs. These datasets were acquired in a previous study for the purpose of simulating intraoperative conditions and analyzing surface deformation \cite{Nakao19}. In this paper, we briefly describe the imaging process and conditions. Contrast-enhanced CT images were measured from the left lungs at two bronchial pressures (14 and 2 cm\ce{H2O}) at the Institute of Laboratory Animals, Kyoto University. This study was performed in accordance with the regulations of the Animal Research Ethics Committee of Kyoto University. All CT images were acquired with a 16-row multidetector CT scanner (Alexion 16; Toshiba Medical Systems, Tochigi, Japan). The dogs were administered ketamine, xylazine, and rocuronium before undergoing tracheal intubation and mechanical ventilation with a ventilator (Savina 300; Drager AG \& Co. KGaA, L\"{u}Beck, Germany). Using the ventilator, bronchial pressure was set to 14 cm\ce{H2O} to obtain images of fully expanded lungs (inflated state) and to 2 cm\ce{H2O} for imaging of collapsed lungs (deflated state). All dogs were placed in a right lateral (decubitus) position on the bed of the CT scanner, and CT images were acquired in the inflated and deflated states. 

\begin{figure}[t]
  \begin{center}
  \includegraphics[width=88mm]{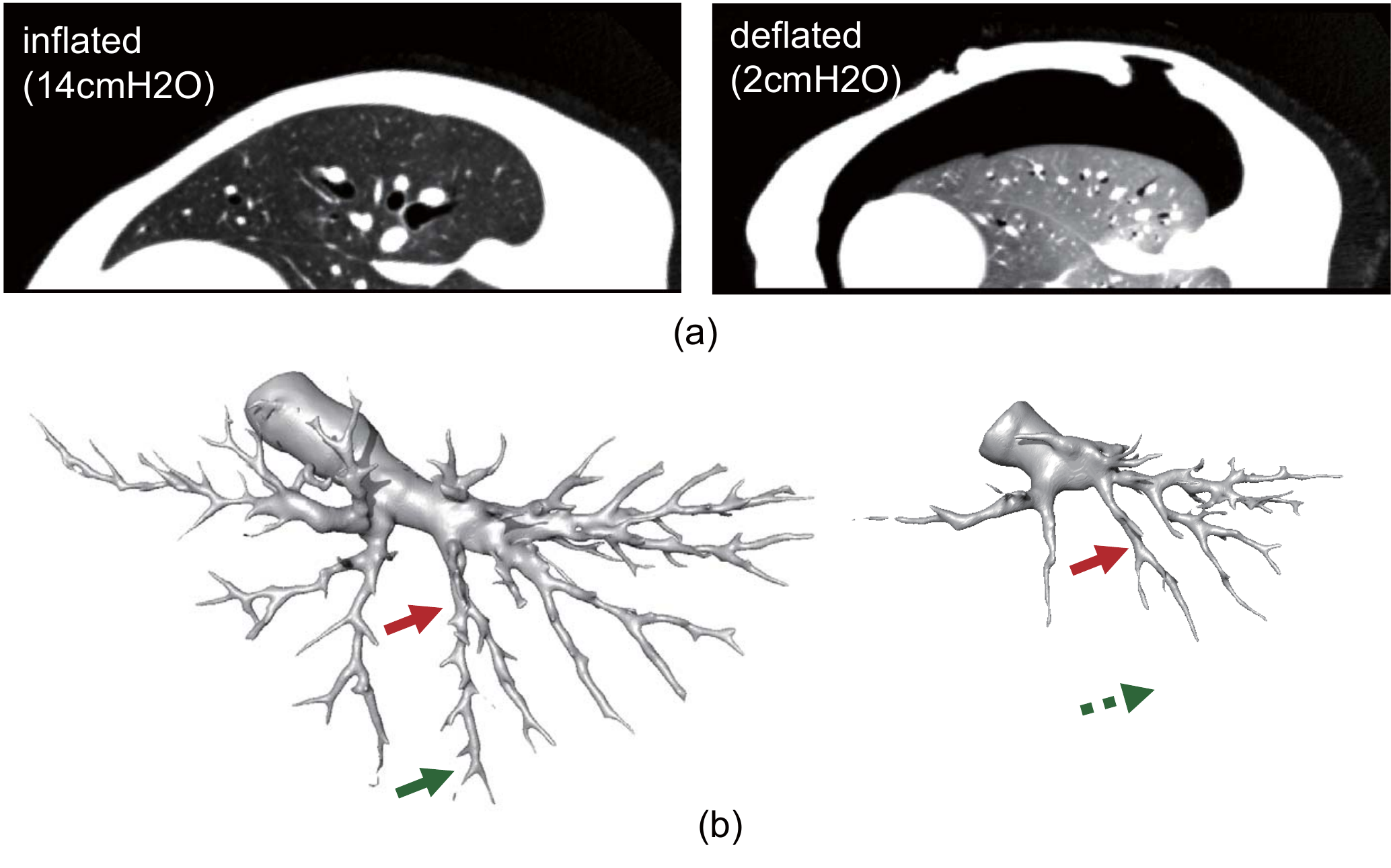}
  \caption{CT slice images and volume rendered results of bronchial structures in the inflated/deflated states. (a) Intensity shift and topological changes inside the thoracic cavity, and (b) the paired bronchial labels extracted using the same CT values. (Arrows indicate estimated corresponding points)} 
  \label{fig:1}
  \end{center}
\end{figure}

% Image characteristics
Fig. \ref{fig:1} shows CT slice images and volume rendered results of bronchial structures in the inflated and deflated states after rigid registration of the two volumes using the spine and bed as fixed references. The window width and level were set to $700$ HU and $-600$ HU, respectively. The CT intensity values changed because of differences in the air content of the lung, with the image contrast in the deflated state becoming significantly lower than that in the inflated state. With this intensity shift, the CT values of the parenchyma increase (the parenchyma region becomes brighter), as can be seen by comparing Fig. \ref{fig:1} (a) and (b). A large space outside of the lung because of air flowing into the thoracic cavity in the deflated state is also confirmed. The accuracy of image-based registration can be affected by these image characteristics, and we found that some image-based registration algorithms \cite{Klein10, Faisal05} failed, having large registration errors in the edge of the deformed lung surfaces. Therefore, analysis of the pneumothorax-associated deformation requires a stable registration framework that addresses the intensity shifts and major deformations with topological changes inside the thoracic cavity. 

\begin{figure*}[t]
  \begin{center}
  \includegraphics[width=165mm]{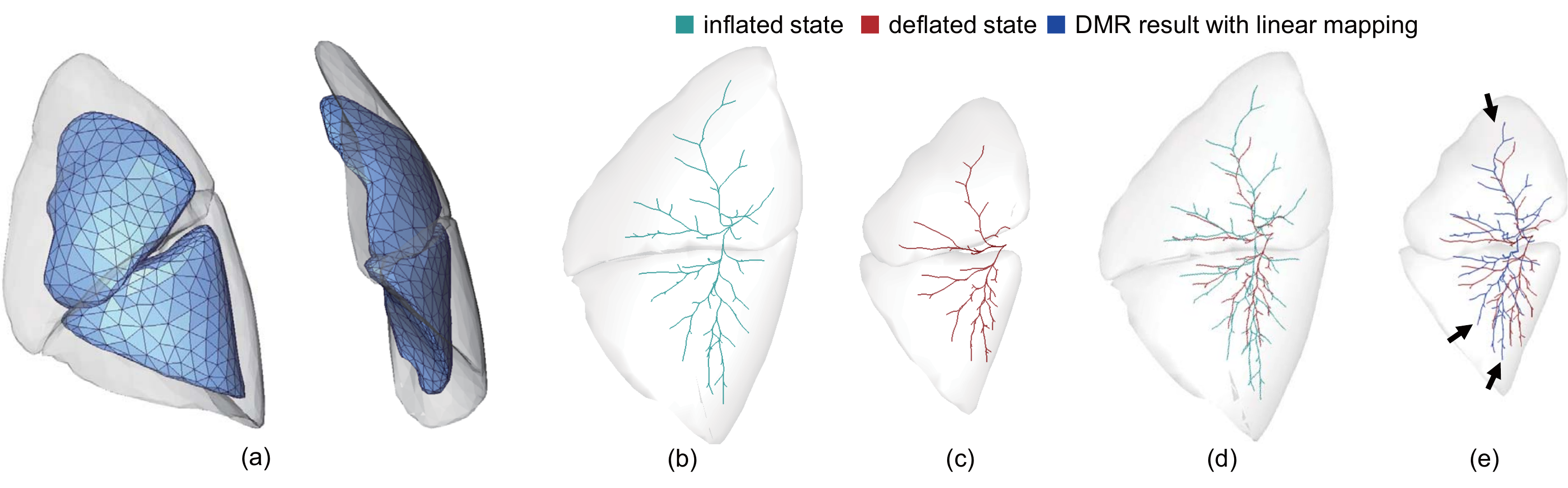}
  \caption{A typical example of the mesh models and DMR results. (a) Tetrahedral meshes in the inflated (translucent) and deflated states (opaque), (b) bronchial centerlines in the inflated state, (c) bronchial centerlines in the deflated state, (d) overlapping images and (e) deformed centerlines (blue) using conventional DMR and linear mapping.}
  \label{fig:2}
  \end{center}
\end{figure*}

\begin{table*}[ht]
	\caption{Specifications of the paired lung geometry datasets in inflated/deflated states. $V_{u}$: volume of upper lobe, $V_{l}$: volume of lower lobe, $N_{j}$: number of bronchial junctions, and $N_{e}$: number of bronchial terminals.}
	\label{table:dataset}
	\centering
	\begin{tabular}{cccccccccc}
		\hline
		\multirow{2}{*}{Case} &\multicolumn{4}{c}{inflated state}& &\multicolumn{4}{c}{deflated state} \\ \cline{2-5} \cline{7-10}
		& $V_{u}$ [cc] & $V_{l}$ [cc] & $N_{j}$ & $N_{e}$ & & $V_{u}$ [cc] & $V_{l}$ [cc] & $N_{j}$ & $N_{e}$ \\
		\hline \hline
		1 & 173.5 & 242.9 & 58 & 62 & & 82.7 & 96.9 & 44 & 48 \\
		2 & 167.1 & 508.3 & 63 & 64 & & 69.9 & 116.4 & 62 & 63 \\
		3 & 231.7 & 357.6 & 66 & 67 & & 111.4 & 145.3 & 40 & 48 \\
		4 & 184.7 & 258.2 & 73 & 75 & & 71.4 & 100.4 & 36 & 40 \\
		5 & 246.1 & 353.3 & 62 & 67 & & 78.9 & 92.8 & 49 & 57 \\
		6 & 188.4 & 373.2 & 77 & 83 & & 65.8 & 104.1 & 30 & 37 \\
		7 & 131.9 & 257.9 & 97 & 98 & & 55.8 & 110.7 & 31 & 32 \\
		8 & 164.3 & 342.1 & 69 & 75 & & 102.3 & 196.7 & 39 & 42 \\
		9 & 198.8 & 375.4 & 58 & 67 & & 84.9 & 111.8 & 43 & 51 \\
		10 & 232.5 & 355.7 & 90 & 98 & & 113.5 & 186.4 & 60 & 64 \\
		11 & 163.7 & 212.1 & $-$ & $-$ & & 160.2 & 236.6 & $-$ & $-$ \\
		\hline
		mean$\pm$SD & 189.3$\pm$35.1 & 330.6$\pm$83.4 & 71.3$\pm$13.3 & 75.6$\pm$13.3 & &  90.6$\pm$29.7 & 136.2$\pm$48.7 & 43.4$\pm$10.9 & 48.2$\pm$10.8 \\
		\hline
	\end{tabular}
\end{table*}

Fig. \ref{fig:1} (c) and (d) show the segmented bronchial labels in the inflated/deflated states, with red and green arrows showing the estimated positions of corresponding points. Although both datasets were extracted using the same CT value threshold, there are clearly missing shapes around the bronchial terminals in Fig. \ref{fig:1}(d). As the inflated state shows the clear contrast of the bronchus because of the high air content, the finer structures can be extracted. This means that deformable registration based on only bronchial structures results in incorrect matching results. We therefore designed our model-based registration methods on the basis of surface and bronchial geometry. 

% Segmentation and modeling
Anatomical segmentation of the upper/lower lobes and bronchial structures was automatically performed using a Synapse VINCENT 3D medical image analysis system (Fujifilm, Tokyo, Japan). Tetrahedral meshes representing the lobes were independently created using Poisson surface reconstruction \cite{Kazhdan06} and tetrahedralization (See Figure 2). All models were resampled to 500 triangles while considering the surface representation and computational cost for model-based registration, resulting in 650--750 tetrahedra for each lobe. The bronchial labels were converted to centerlines using the Vascular Modeling Toolkit (VMTK) \cite{VMTK}. Figure 2 (b) and (c) show the bronchial centerlines of inflated and deflated states, respectively. The centerlines comprised a tree structure with multiple nodes, including junctions and terminals. The root node represents the pulmonary hilum. TABLE \ref{table:dataset} summarizes the specifications of the 11 paired lung models in the inflated and deflated states. $V_u$ and $V_l$ are the volumes of the upper and lower lobes, respectively. $N_j$ is the number of bronchial junctions, and $N_e$ is the number of bronchial terminals. These numbers are determined in VMTK based on the complexity of the centerline structures, with the deflated models showing smaller values. We note that in Case 11, the volume at a bronchial pressure of 2 cm\ce{H2O} showed an unexpected increase, which was probably a result of the bronchial pressure not being successfully controlled. Therefore, Case 11 was removed from our analyses.

\subsection{Problem definition}

% Model-based registration
The focus of this study was to investigate the spatial displacements of whole lungs, including their surfaces and bronchial structures. The displacement of bronchi can be partially obtained by matching anatomical feature points extracted from bronchial junctions \cite{Nakamoto07, Uneri13}. However, such landmark-based matching cannot be applied to the lung’s curved surfaces because of their lack of anatomical features. Model-based registration and DMR are increasingly being applied in practice because registered mesh models are directly available for statistical modeling and variational analysis \cite{Ehrhardt11, Kim15, Jud17}. In addition, for anatomical structures with large shape variations, recent DMR methods \cite{Kim15, Nakao19} have shown better registration accuracy than large-deformation diffeomorphic metric mapping \cite{Faisal05}.

\begin{figure*}[t]
	\begin{center}
		\includegraphics[width=175mm]{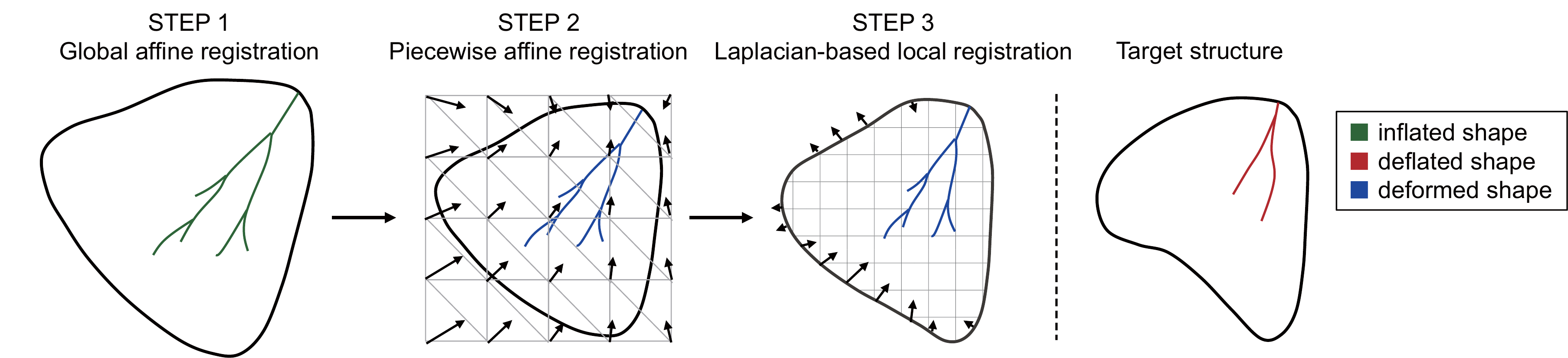}
	\end{center}
	\caption{Proposed three-step DMR framework for mixed data structures with lung surfaces and centerlines. The registration process is applied to each lobe independently to cope with rotational components or sliding motion of pneumothorax-associated deformation.}
    \label{fig:3}
\end{figure*}

% Deformation analysis through DMR
The pair of lung mesh models representing the inflated and deflated states differed in the number of vertices and the data structure (i.e., tree and mesh topology) because they were independently generated from different CT images. In this study, the shape model $M_I$ generated from the inflated lung is deformed to the deflated model $M_D$. We aimed to compute the registered model $\phi(M_I) $ ($\phi$: mapping function) that precisely approximated the target surfaces and centerlines by DMR. Because the initial and registered model $(M_I, \phi(M_I))$ have the same data structure with point-to-point correspondence, the spatial deformation $\Vec{u}$ can be obtained by calculating the displacement vectors of the corresponding vertices. Unlike organ segmentation, which aims to represent the anatomical region, accuracy of local correspondence is required at each vertex of the mesh. 

% Heterogeneity
We performed a preliminary experiment to investigate the linearity of displacement in the lung lobes. We applied conventional DMR\cite{Nakao19} to the two lung surfaces of the inflated/deflated states and mapped the vertex displacement of the deformed tetrahedral mesh to all the nodes of the bronchial tree structures. This mapping involved a combination of DMR and spatial linear interpolation of the inside of the lobe surfaces and was implemented using the barycentric coordinates of the tetrahedral elements. Fig. \ref{fig:2} (e) shows the DMR result with linearly-mapped bronchial centerlines (blue) and the target centerlines  in the deflated state (red). 

% Problem definition
If the pneumothorax-associated deformation is spatially homogeneous, the deformed bronchial centerline would coincide with the target centerline in the deflated state. However, as shown in Fig. \ref{fig:2}(e), local mismatches of the bronchial junctions were found, as indicated by the arrows. Most of the datasets showed large deviations between the target and linearly mapped centerlines, suggesting that the displacement of the lung region was not spatially homogeneous. This preliminary analysis also shows that the pneumothorax-associated deformation is spatially continuous in each lobe, although some cases show irregular rotations between the upper and lower lobes. Based on the features observed between the lung surface and the bronchi, we designed the DMR framework by considering the following two issues.

\begin{itemize}
    \item The lung surface and bronchial centerlines should be registered under the constraint that the deformation is spatially continuous and smooth while matching visible junctions as much as possible. 
    \item To handle the different rotations between the upper and lower lobes, mapping functions $\phi$ for deformable registration should be determined independently for each lobe. 
\end{itemize}

\subsection{DMR for lung surface and bronchial structures}
Based on the problem definition, this section describes the outline of the DMR framework that matches the inflated and deflated lung models. The proposed DMR is designed to obtain non-linear mapping under the condition that the deformation field for the lung surface and bronchi is spatially continuous and smooth. The registration framework comprises the following three-step processes (See Figure \ref{fig:3}).

\begin{enumerate}
    \item [] \hspace{-7mm} STEP 1 Global affine transformation \\
    Global registration is first applied to the inflated (source) model $M_I = (S_I, C_I)$ ($S_I$: the lung surface and $C_I$: the bronchial centerline in the inflated state) to roughly match it to the deflated (target) model $M_D = (S_D, C_D)$ ($S_D$: the lung surface and $C_D$: the bronchial centerline in the deflated state). In this step, a set of vertices of the lung surface and bronchial centerlines are transformed using $ \Vec{v}_i' = T\Vec{v}_i$ ($T$: transformation matrix, $\Vec{v}_i$: vertex position). As the root node of the bronchial tree structure represents the pulmonary hilum in our data, its position is used as the origin for normalizing the translation.

    \item [] \hspace{-7mm} STEP 2 Piecewise affine (PWA) deformable registration \\
    The purpose of the next step is to locally register the surface and bronchial centerline while spatially deforming the entire lung region continuously and smoothly. A cuboid grid $\mathcal{M}$ formed by tetrahedral elements is first generated as the bounding box of the lung model $\phi(M_I)$ obtained in STEP 1. The displacement of each vertex of the grid is mapped to all the surface and centerline vertices, representing the deformation fields. We introduce an objective function that addresses local mismatches on the missing parts of the bronchi while handling the mixed data structures with surface meshes and centerlines. The details of the objective function are described in the next section.
    
    \item [] \hspace{-7mm} STEP 3 Laplacian-based local mesh registration \\
    The final step is performed to improve registration results for the local shapes with large matching errors. The lung surface is refined based on the Laplacian-based diffeomorphic shape matching (LDSM) algorithm \cite{Nakao20}, and the deformation is mapped to the centerlines. Using a discrete Laplacian of the lung surface as the shape descriptor, LDSM achieves stricter surface matching while preserving the entire shape of the lung. Details of the LDSM algorithm and its registration performance for abdominal organs are reported in \cite{Nakao20}. 
\end{enumerate}

We applied the above registration process independently to each lobe to cope with rotation and sliding motion of the upper and lower lobes. Hence, the transformation matrix $T$ for global affine transformation and mapping functions $\phi$ (i.e., deformation field $\Vec{u}$) for deformable registration was calculated for each lobe, and simultaneously updated per iteration in the registration process. For the propagation of deformation fields to the surface and bronchial structures, the barycentric coordinates of the vertices with regard to the corresponding tetrahedral element were calculated in advance. All vertices were updated while keeping the relative position of the four points of the tetrahedral element. 

\subsection{Objective function}
In STEP 1 and 2, the registered models are computed by minimizing the proposed objective function, which is designed to obtain a spatially continuous and smooth mapping of the lung surface and bronchial centerlines. The objective function is described as follows.

\begin{align}
    E &= d(M_D, \phi(M_I)) + \int_{0}^{1} \|L(\Vec{u}(s))\|^2 ds,
\end{align}
where $d$ is the distance function between the two models. $\phi(\cdot)$ is a continuous and differentiable transformation that maps the source model to the deformed model. $L(\cdot)$ is the Laplace--Beltrami operator and $L(\Vec{u})$ is the discrete Laplacian of the deformation field \Vec{u}.

\begin{figure}[t]
	\begin{center}
		\includegraphics[width=84mm]{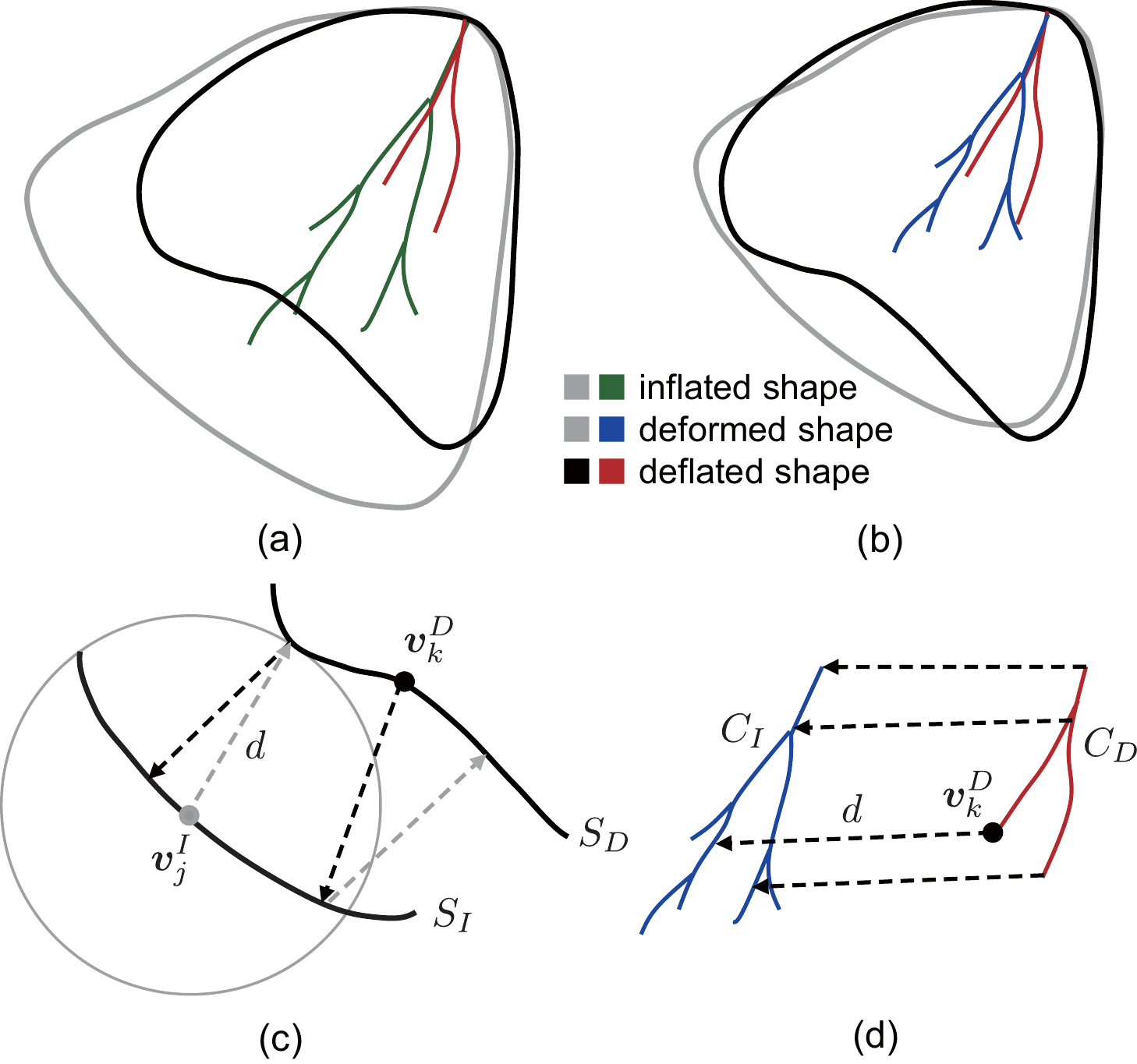}
	\end{center}
	\caption{Local distance measures for lung surfaces and bronchial centerlines with missing parts. (a) initial setup of the inflated and deflated models, (b) deformed and deflated models after global affine transformation, (c) bidirectional surface distance, and (d) one-way centerline distance.}
    \label{Fig:4}
\end{figure}

The first term evaluates the difference between the deformed model and the target model, whereas the second is a regularization term to make the deformation field smooth. Figs \ref{Fig:4} (a) and (b) illustrate examples of initial setups in Steps 1 and 2, respectively. To obtain accurate deformation fields that match the pneumothorax lung models with partial defects, we introduce the following distance measures.

\subsubsection{Surface distance} 
The surface distance $d_s$ is a loss function that evaluates the difference between the two curved surfaces. In this study, the distance between the source and target surfaces $(S_I, S_D)$ is defined as follows:

\begin{align}
    \label{eq:loss1}
    d_{s}(S_I, S_D) % &= \int d(S_I(s), S_D) ds + \int d(S_D(s), S_I) ds \nonumber \\
    &= \sum_{\Vec{v}_j^I \in S_I} d(\Vec{v}_j^I, S_D) / N_I + \sum_{\Vec{v}_k^D \in S_D} d(\Vec{v}_k^D, S_I) / N_D,
\end{align}
where $S_I$ and $S_D$ are the set of vertices for the inflated and deflated lung surfaces, respectively, and $N_I$ and $N_D$ are their number of vertices. Although the source surface is iteratively updated in the registration process, we use $S_I$ for its notation, for readability purposes. $\Vec{v}_j^I$ represents the 3D position of the $j$-th vertex of $S_I$, and $\Vec{v}_j^I$ is the 3D position of the $k$-th vertex of $S_D$. $m$ and $n$ are the numbers of vertices for $S_I$ and $S_D$, respectively. $d(\Vec{v}, V)$ is the minimum point distance between a vertex ${\Vec{v}}$ and a set of vertices ${V}$ of a curved surface, and is defined by

\begin{align}
    \label{eq:vs_dist}
    d(\Vec{v}, S) &= | \Vec{v} - \hat{\Vec{v}_p} | \nonumber \\ 
    s.t. \hspace{3mm} \hat{\Vec{v}_p} &= \min_{\Vec{v}_p} \left\{  | \Vec{v} - \Vec{v}_p |  + \gamma ( 1 - \Vec{n} \cdot \Vec{n}_p ) \right\},
\end{align}
where $\gamma$ is a weight, $\Vec{v}_p$ is a vertex of the curved surfaces $S$, and $\Vec{n}$ and $\Vec{n}_p$ are vertex normals at $\Vec{v}$ and $\Vec{v}_p$, respectively. To avoid matching errors around the ridges or tips of the lung surface, the closest point $\hat{\Vec{v}_p}$ from $\Vec{v}$ is determined by considering the vertex normal. Figure 4(c) illustrates $d(\Vec{v}_j^I, S_D)$ (gray) and $d(\Vec{v}_k^D, S_I)$ (black). $d_s$ is computed as the average of the bidirectional distances from the surface to the other surface and is available as a stable index even for highly frequent  shape changes and noise. This error metric is also referred to as a mean distance and is used to evaluate registration accuracy\cite{Kim15, Nakao19}. \\

\subsubsection{One-way centerline distance} 
The one-way centerline distance $d_c$ is a loss function that evaluates the difference between the two centerlines $(C_I, C_D)$ of the source and the target model. Unlike for the surface distance, we employ a one-way distance measure to consider the missing shapes of the deflated bronchi. The distance is described by

\begin{align}
    \label{eq:loss2}
    d_{c}(C_I, C_D) % &= \int d(C_I(s), S_D) ds
    &= \sum_{\Vec{v}_k^D \in C_D} d(\Vec{v}_k^D, C_I) / N_C,
\end{align}
where $C_I$ and $C_D$ are the set of vertices for the bronchial centerlines in the inflated and deflated states, $N_C$ is the number of vertices for $C_D$, $d(\Vec{v}, C)$ is the minimum point distance between a vertex ${\Vec{v}}$ and a curved line ${C}$. Figure 4(d) shows $d(\Vec{v}_k^D, C_I)$ uniquely defined as the one-way distance from a point of the target centerline to the source point. The source centerline $C_I$ has many terminals that do not structurally correspond to those of the target centerline $C_D$ because of its missing branches. Therefore, the bidirectional distance between the two centerlines results in an incorrect evaluation. We can assume that the source centerline $C_I$ has a more complex structure than the target centerline because of the CT intensity shift between the inflated and deflated states. The proposed loss function enables evaluation of the validity of the  distances corresponding to the visualized junctions and branches. \\

\subsubsection{Full objective}
The discrete form of the objective function $E$ can be described using the two distance measures as follows:

\begin{align}
    E = d_{s}(S_I, S_D)^2 + \alpha d_{c}(C_I,C_D)^2 + \beta \sum_{\Vec{v}_i \in M_I} \|L(\Vec{u}_{i})\|^{2},
    \label{eq:full_obj}
\end{align}
where $\alpha$ and $\beta$ are the weights, $\Vec{u}_{i}$ is the local displacement at vertex ${\Vec v_i}$, and $L({\Vec u_i})$ is the discrete Laplacian of the local displacement vector, described as

\begin{align}
	L(\Vec u_{i}) = \sum_{\Vec{v}_j \in A_{i}} \omega_{ij} (\Vec u_{i}- \Vec u_{j}).
\end{align}
 Here, $\omega_{ij}$ is the edge weight and $A_{i}$ is the set of adjacent vertices of one ring connected by vertex ${\Vec v_i}$ and the edges. As the discrete Laplacian approximates the mean curvature normal\cite{Nealen06}, the third term in Eq. (\ref{eq:full_obj}) makes the deformation field spatially smooth and diffeomorphic. Although there are several methods for calculating the weights, the general one is cotangent discretization based on the per-edge Voronoi areas. The deformation field $\Vec{u}$ is obtained by iteratively updating $\Vec{v}_{i}$ while minimizing $E$. The step-by-step update avoids local mismatches during the early stage if there is a considerable distance between the two models.

\subsection{Contraction and rotational components}
This section describes methods for analyzing the heterogeneity of the pneumothorax-associated deformation using the registered models. Figure 5(a) illustrates an ideal case of homogeneous deformation with zero rotational components. As the pulmonary hilum is physically fixed and not affected by pneumothorax-associated deformation, we use its position as the origin. In this case, the per-vertex displacement $\Vec{u}_i = \phi(\Vec{v}_i
)
 - \Vec{v}_i$ towards the hilum corresponds to the local contraction component $\Vec{s}_i$, and we can estimate the local characteristics of the surface contraction using the size of $\Vec{u}_{i}$.
 
Figure 5(b) shows a real example of pneumothorax-associated deformation with rotation. The contraction was accompanied by various types of rotations, as observed in different subjects. As each subject was measured in the same right lateral position, we consider it unlikely that posture differences were the major cause of rotational variations.  Rather, we consider that such rotations are caused by sliding motions between lobes and subject-specific physical interactions between the lobes and the thoracic cavity, such as friction and adhesion. Therefore, the rotational components can be regarded as external disturbance that does not reflect differences in the characteristics of lung tissue between subjects. As this study aimed to analyze the tissue characteristics of the lung lobes, rotational components caused by external factors were excluded from our analysis. 

\begin{figure}[t]
	\begin{center}
		\includegraphics[width=80mm]{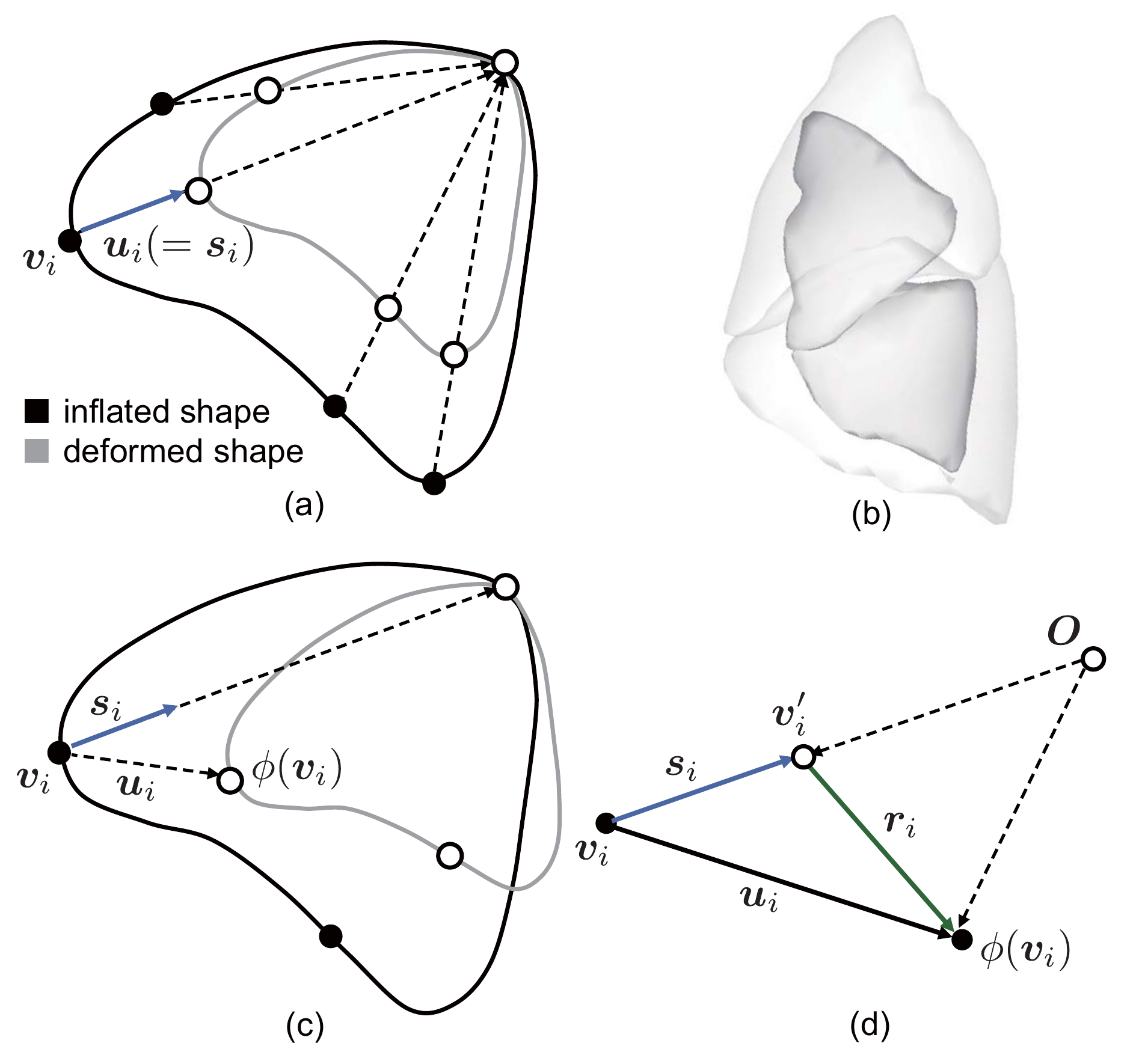}
	\end{center}
	\caption{Contraction and rotational components in pneumothorax-associated deformation. (a) Homogeneous deformation with zero rotational components, (b) a real example of pneumothorax-associated deformation with rotation, (c) the contraction component extracted from the deformation field, and (d) its geometrical description.}
	\label{fig:5}
\end{figure}

We employed simple calculation methods to extract the contraction component from the deformation field. In Fig. 5 (c), $\Vec{v}_i$ and $\phi(\Vec{v}_i)$ are corresponding vertices of the initial model and the registered model determined by DMR, respectively. To extract contraction $\Vec{s}_i$ from the displacement vector $\Vec{u}_i$, $\Vec{v}'_i$ is defined as the position where the rotational component $\Vec{r}_i$ is removed from $\phi(\Vec{v}_i)$. As the pulmonary hilum (i.e., the origin in this study) is supposed to be the center of rotation, the direction of the contraction component $\Vec{s}_i$ is defined as the vector towards the hilum, as shown in Figure \ref{fig:5}(d). The contraction $\Vec{s}_i$ and rotational component $\Vec{r}_i$ are defined as Eq. (\ref{eq:contraction}) and (\ref{eq:rotation}), respectively. 

\begin{align}
    \Vec{s}_i=\frac{|\phi(\Vec{v}_i)|-|\Vec{v}_i|}{|\Vec{v}_i|} \Vec{v}_i       
    \label{eq:contraction},
\end{align}

\begin{align}
    \Vec{r}_i=\Vec{u}_i-\Vec{s}_i.
    \label{eq:rotation}
\end{align}

\subsection{Strain definition for bronchi and parenchyma}
The objective of this study was to analyze the heterogeneity in pneumothorax-associated deformation of the lungs. For this purpose, the strains of bronchi and other parenchyma regions were compared according to the registered CT models. In general, strain is a measure of the displacement between particles in the body relative to a reference length. The Cauchy strain is defined as  

\begin{align}
    \varepsilon= \frac{\Delta L}{L}
    \label{eq:strain}
\end{align}
where $L$ is the original length between two sampled points, and $\Delta L$ is the change in length between the sampled points in the deformed body. To compare the strains of anatomically corresponding structures, we divided the lung's subsurface structures into two regions: the bronchial region and the region between the bronchial terminal and the surface (notated as parenchyma in this paper). This latter region comprises microtissues such as alveoli and terminal bronchioles, and the structures are not clearly identifiable on CT imaging. To determine $L$ and $\Delta L$ for the corresponding two regions, we focused on a set of spatially-continuous points per bronchial branch, bronchial terminal, and its extended point on the lung surface. Figure 6(a) shows the sampled points: the pulmonary hilum, three bronchial junctions (blue), the bronchial terminal (green), and corresponding surface point (red) obtained from one bronchial branch. The surface point is defined as the intersection of the straight line determined by the pulmonary hilum and the bronchial terminal. In this case, $L$ is determined as the Euclidean distances of the two neighboring points in the deflated state and $\Delta L$ is the change in the distance that occurs in the deformation between the inflated state and the deflated state. 

Figure 6 (b) illustrates a 2D plot representing the relation between the Euclidean distance between the pulmonary hilum $|\Vec{v}_i|$ and the contraction $|\Vec{s}_i|$ for the sampled points. Here, the gradient of the line defined by the two neighboring points represents the Cauchy strain for the local region. As we confirmed that the set of plots for one bronchial branch were mostly distributed linearly, the strain for the bronchi was calculated from the linear regression of the plots from the bronchial junctions and the terminal, as shown by the blue line in Fig. 6(b). However, the strain for the region between the bronchial terminal and the lung surface was obtained from the gradient of the red line.

\begin{figure}[t]
	\begin{center}
		\includegraphics[width=80mm]{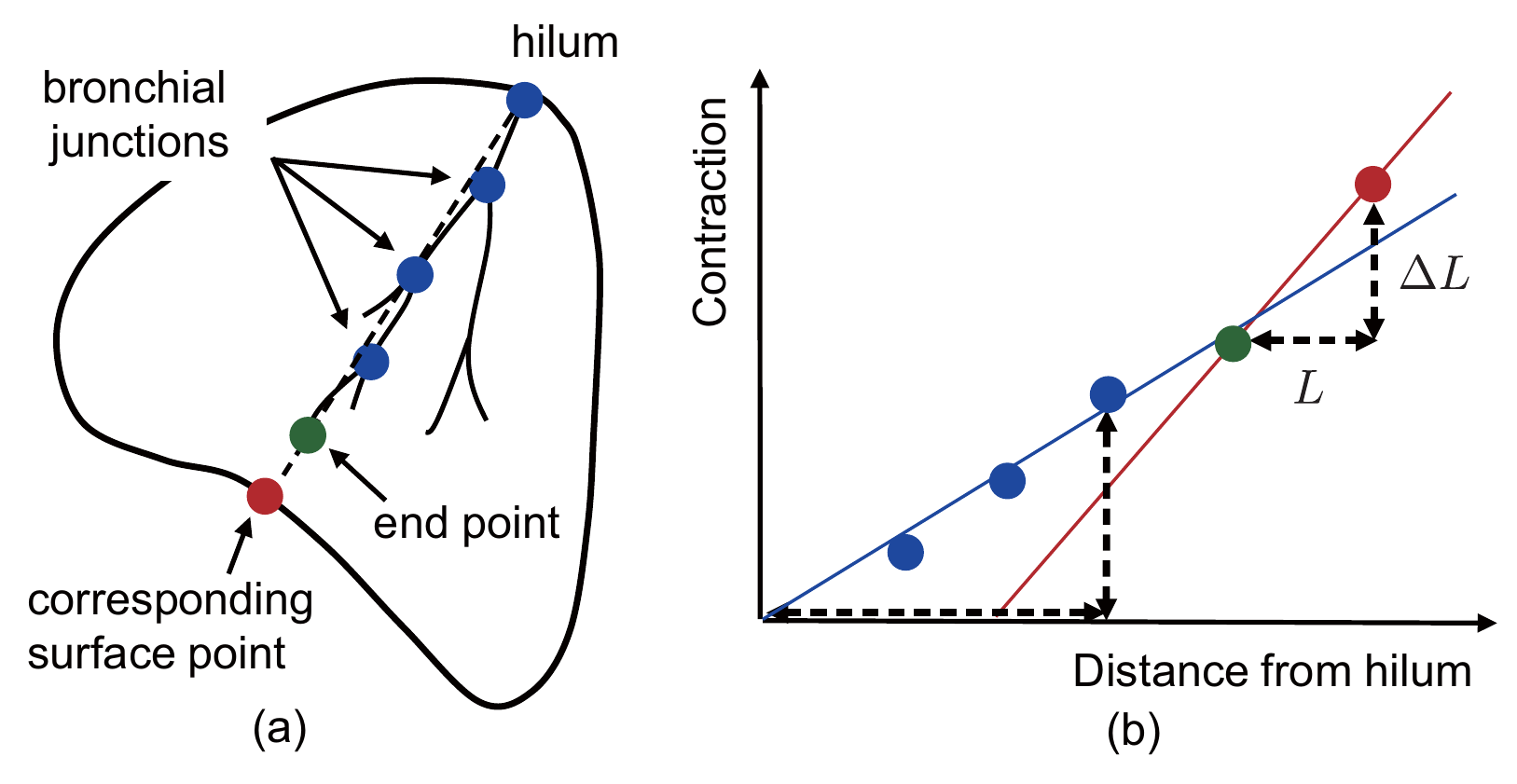}
	\end{center}
	\caption{Definition of the Cauchy strain for bronchi and parenchyma in each lobe. (a) The positional relationship between the hilum, bronchial junctions, and terminal and corresponding surface point. (b) 2D plots and strain calculation for bronchi (blue line) and parenchyma (red line).}
\end{figure}

\section{Experiments}
In the experiments, the performance of the proposed DMR framework was first evaluated using ten in vivo lung models (Cases 1--10). The efficacy of the objective function was confirmed by comparing it with existing model-based registration approaches. Then, the heterogeneity of the strains on the bronchus and parenchyma was investigated using the spatial displacement obtained from the registered models. The proposed DMR and deformation analysis methods were implemented using Visual C/C++ and OpenGL. A computer with a graphics processing unit (CPU: Intel Core i7 3.7GHz, Memory: 64 GB, GPU: NVIDIA GeForce GTX 1080) was used throughout the experiments. 

% parameter and convergence settings
For the weights in the developed framework, we used 1.0 for $\lambda$, the same value used in a previous study \cite{Nakao19}. After examining several parameter sets, $(\alpha, \beta) = (2.0, 2.0)$ were used for the weight parameters in the objective function. For each registration step, when the minimum value of the objective function was not updated in the most recent 20 iterations, or when the number of updates reached 1000 times, the iterative update was terminated and the next registration step was initiated. 

\subsection{Quantitative comparison of DMR}

% Shape similarity metrics
In the first experiment, the registration accuracy was confirmed for the paired CT models with pneumothorax-associated deformation. The inflated lung models were deformed to the deflated target models and the differences between the deformed and target shapes were quantified using multiple shape similarity criteria. The mean distance and Hausdorff distance \cite{Kim15, Nakao19} were used as measures of the shape similarity of the surfaces. The mean and maximum of the one-way point-to-point distances were used to measure the similarity of the bronchial centerlines. The Hausdorff distance measures the longest distance among the minimum point distances between two surfaces, whereas the mean distance is the average of the minimum point distances. The Hausdorff distance is a strict error metric as it measures the largest distance between two shapes. 

\begin{figure*}[t]
  \begin{center}
  \includegraphics[width=180mm]{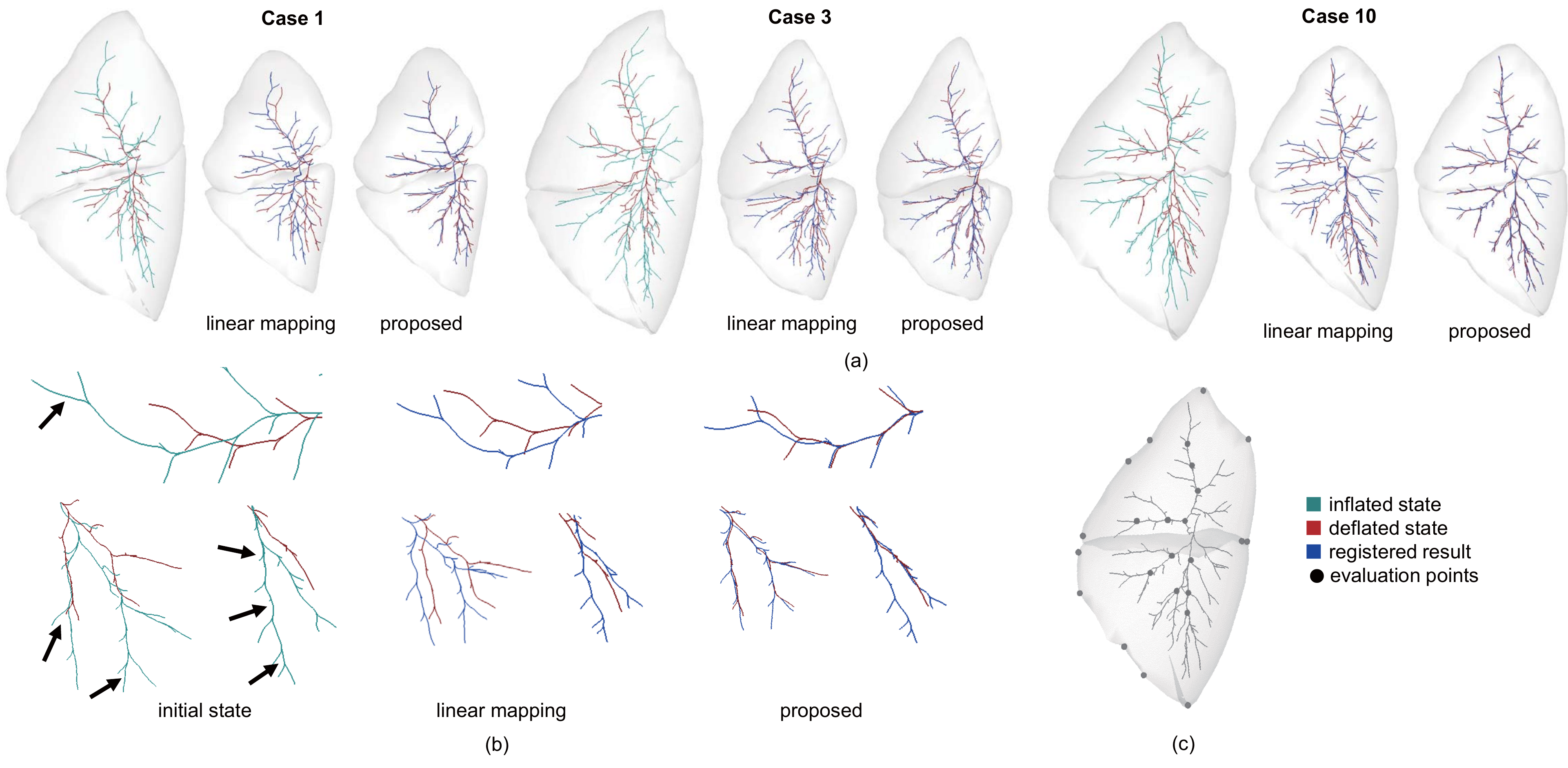}
  \caption{Typical registration results with variations in bronchial anatomy and lobe deformation. (a) visual comparison between the proposed DMR and surface registration with linear mapping, (b) enlarged images of the local bronchial shapes, (c) 24 points for evaluating target registration error.}
  \label{fig:registration_results}
  \end{center}
\end{figure*}

\begin{table*}[t]
	\caption{The mean $\pm$ standard deviation of DMR results for mean distance (MD), Hausdorff distance (HD), centerline distance (CD), and target registration errors (TREs) for surfaces and bronchial centerlines.}
	\label{table:registration}
	\centering
	\begin{tabular}{cccccc}
		\hline
	    Error metric & Rigid &  AF & PWA & LSM & Proposed \\ 
		\hline \hline
		MD [mm] & 10.09 $\pm$ 2.14 & 4.13 $\pm$ 1.22 & 1.49 $\pm$ 0.20 & 0.24 $\pm$ 0.05 & 0.24 $\pm$ 0.06 \\
		HD [mm] & 29.65 $\pm$ 5.72 & 15.97 $\pm$ 4.84 & 6.35 $\pm$ 1.85 & 0.74 $\pm$ 0.18 & 0.78 $\pm$ 0.17 \\
		CD (mean) [mm] & 4.08 $\pm$ 1.15 & 3.83 $\pm$ 1.80 & 3.88 $\pm$ 1.68 & 3.86 $\pm$ 1.70 & 1.51 $\pm$ 0.58* \\
		CD (maximum) [mm] & 11.48 $\pm$ 4.71 & 12.78 $\pm$ 6.61 & 12.07 $\pm$ 6.04 & 12.40 $\pm$ 6.14 & 7.53 $\pm$ 2.47* \\
		TRE (surface) [mm] & 22.91 $\pm$ 8.81 & 10.65 $\pm$ 5.46 & 6.11 $\pm$ 3.69 & 4.57 $\pm$ 4.13 & 4.55 $\pm$ 2.82 \\
		TRE (bronchus) [mm] & 7.59 $\pm$ 4.77 & 5.63 $\pm$ 3.60 & 6.10 $\pm$ 3.54 & 6.14 $\pm$ 3.51 & 3.97 $\pm$ 2.75* \\
		\hline
		& & & & & \multicolumn{1}{r}{*: $p<$0.01}\\
	\end{tabular}

\end{table*}

% Target registration error (TRE)
Unlike segmentation or recognition problems, deformation analysis requires point-to-point correspondence between two shapes. The accuracy of the registered shape after DMR was validated by the target registration errors (TREs) on anatomical evaluation points of the lung surface and bronchial junctions. The TREs were defined as the Euclidean distances between the manually-placed evaluation points of the target (deflated) model and the deformed (registered) model generated from the inflated model. Automatic identification of the unique/distinctive features of individual models is difficult because their shapes comprise smooth curved surfaces and their bronchial tree structures vary substantially among patients (See Figure \ref{fig:registration_results}). With the anatomical review of two expert surgeons, 12 evaluation points were manually set on the surfaces and bronchial junctions of the deflated models. The corresponding vertices of the inflated models were also determined. Figure \ref{fig:registration_results}(c) shows an example of all 24 evaluation points for one subject. The evaluation points for the surfaces were selected from the ridge of the lung surface with large curvatures and locally similar features. As the bronchial tree structures were anatomically different between subjects, as shown in Table \ref{table:dataset}, the bronchial junctions were selected such that the correspondence between the inflated/deflated shapes was visually clear.

% Reference methods
The proposed DMR was compared with the existing model-based registration methods:

\begin{itemize}
    \item Affine transformation (AF)
    \item Piecewise affine transformation (PWA) 
    \item Laplacian-based shape matching (LSM) 
\end{itemize}

The rigid registration result was also listed to confirm the initial difference and complexity of the inter-subject matching of collapsed lungs. PWA and LSM were investigated as an integrated registration framework to match the experimental conditions with those of the proposed method; for PWA and LSM, affine registration was first applied to globally match the posture and volume of the overall shape. Additionally, in the LSM, piecewise affine registration was performed in advance to achieve a globally stable and locally strict registration, similar to the proposed method. As these existing methods were originally applied to surface models, subsurface registration was performed with linear mapping, that is, the displacement of the bronchial structures was linearly interpolated from the deformed surfaces. In the PWA step, an axis-aligned bounding box was calculated for each lobe, and its space was divided into $4 \times 4 \times 4 = 64$ subspace. We note that the difference between the LSM and proposed method lies in the weight parameters $(\alpha, \beta) = (0.0, 2.0)$ used for the objective function in LSM. This setting also serves as confirmation of the role of the proposed one-way centerline distance. 

% Registration performance
Table \ref{table:registration} shows values of the error metrics obtained from the registration results of the 10 subjects. The mean distance (MD), Hausdorff distance (HD) for surfaces, the mean and maximum centerline distance (CD), and target registration errors (TREs) for lung surface and bronchi are listed. The proposed DMR achieved significantly smaller CDs and TREs for bronchial centerlines compared with the conventional linear mapping (one-way analysis of variance, ANOVA; $p < 0.01$ significance level).  There was no significant difference in MD, HD, and TRE for lung surfaces, which means that the proposed DMR outperformed the conventional linear mapping in matching bronchi with missing branches, while maintaining a similar registration quality for surfaces. Specifically, the small TREs of less than 5 mm demonstrate that spatial displacement analysis focusing on the difference between the lung’s bronchi and parenchyma tissues is possible using the registered models.

\begin{figure*}[t]
  \begin{center}
  \includegraphics[width=180mm]{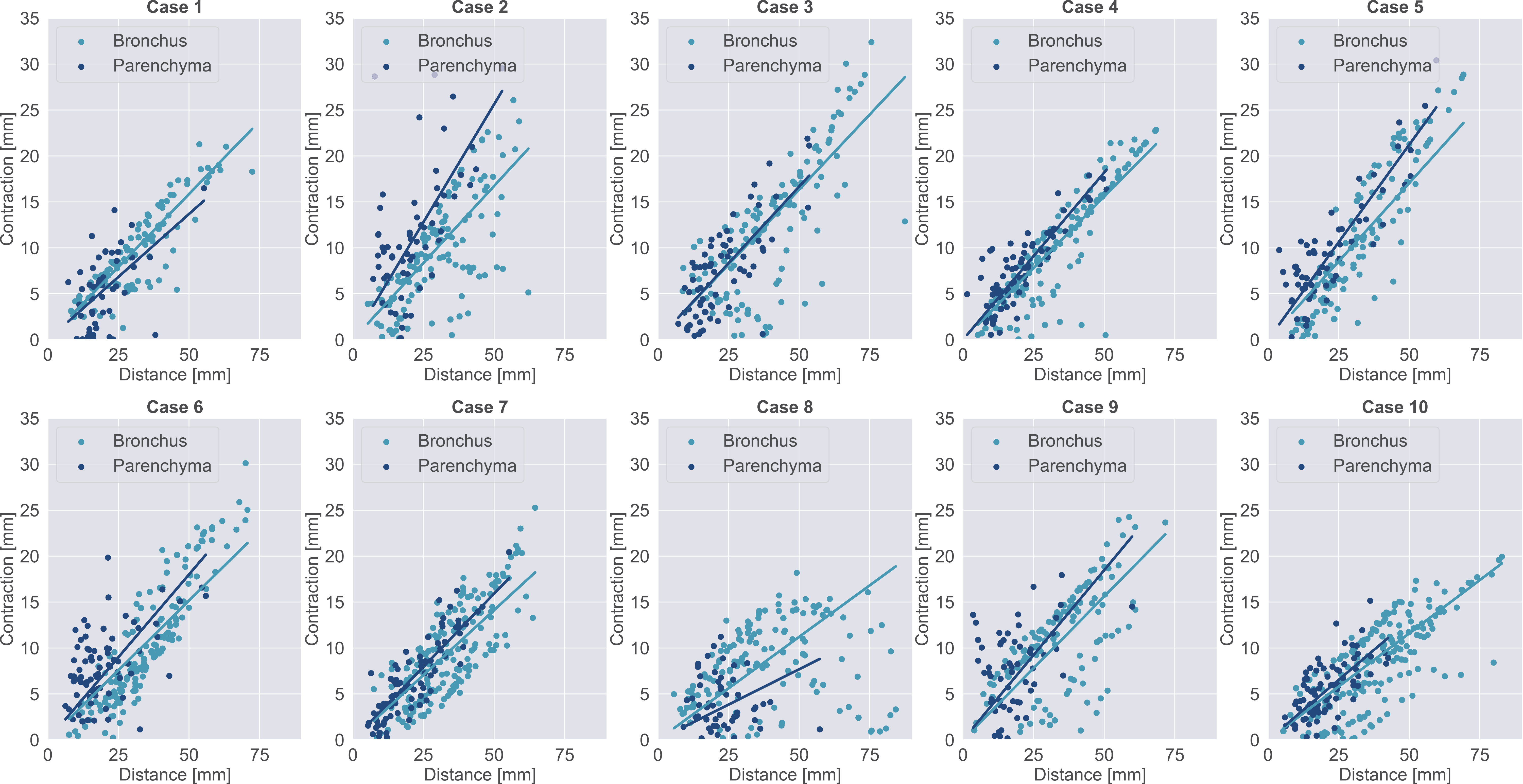}
  \caption{The relationship between the Euclidean distance from the pulmonary hilum and the contraction components for bronchi and parenchyma. The gradient of the regression line represents the strain.}
  \label{fig:strains}
  \end{center}
\end{figure*}

% Registration Examples
Figure \ref{fig:registration_results}(a) shows typical registration results of three subjects (Cases 1, 3, and 10), which include variations in the large deformation with rotations and missing shapes in the bronchial structures. It can be visually confirmed that compared with the conventional method, the proposed DMR achieves smaller errors in matching bronchi. Figure \ref{fig:registration_results} (b) shows enlarged images of the local bronchial shapes for Cases 1 and 6. There are some missing shapes in the target bronchi, as indicated by the arrows. Point-to-point registration using only bronchial tree structures resulted in large registration errors including irregular bends and overfitting due to the underdetermined matching problem. In the proposed method, such errors or mismatching results were not confirmed in all subjects. The proposed objective function evaluates the deviation of both the surface and bronchial structures, and spatially-valid registration is performed, even for the target models with missing information.

\subsection{Analysis of pneumothorax-associated deformation}
The strains of the sampled points of the bronchi and parenchyma were investigated using the analysis method described in Section II-E. As all bronchial branches were targeted, the number of branches was equal to the number of terminals $N_t$ of the deflated model in Table \ref{table:dataset}. For each branch, the displacement vectors for the bronchial junctions, terminals, and corresponding surface points were obtained from the registered models, and the contraction component derived after removing irregular rotation was plotted.

\begin{table}[t]
	\caption{The mean $\pm$ standard deviation of strains for bronchus and parenchyma in pneumothorax-associated deformation}
	\label{table:strain}
	\centering
	\begin{tabular}{ccc}
		\hline
        Case & Bronchus & Parenchyma\\
		\hline \hline
		1 & 0.315 $\pm$ 0.064 & 0.286 $\pm$ 0.194 \\
		2 & 0.372 $\pm$ 0.176 & 0.616 $\pm$ 0.506 \\
		3 & 0.320 $\pm$ 0.209 & 0.329 $\pm$ 0.162 \\
		4 & 0.275 $\pm$ 0.099 & 0.430 $\pm$ 0.398 \\
		5 & 0.247 $\pm$ 0.138 & 0.499 $\pm$ 0.330 \\
		6 & 0.230 $\pm$ 0.109 & 0.463 $\pm$ 0.239 \\
		7 & 0.250 $\pm$ 0.131 & 0.303 $\pm$ 0.138 \\
		8 & 0.295 $\pm$ 0.102 & 0.209 $\pm$ 0.159 \\
		9 & 0.350 $\pm$ 0.128 & 0.580 $\pm$ 0.640 \\
		10 & 0.255 $\pm$ 0.121 & 0.272 $\pm$ 0.118 \\
		\hline
		Mean $\pm$ SD & 0.292 $\pm$ 0.141 & 0.395 $\pm$ 0.296* \\
		\hline
		\hline
		& & \multicolumn{1}{r}{*: $p<$0.05}\\
		\end{tabular}
\end{table}

Figure \ref{fig:strains} is a 2D plot showing the relationship between the reference distance and the contraction components for bronchi and parenchyma. The Euclidean distance from the pulmonary hilum was used for the reference distance of the bronchi, and the Euclidean distance from the bronchial terminal was used for the parenchyma. The number of plots in the graph equals the number of sampled points for all cases. The light and dark blue plots are the contractions of the bronchi and parenchyma, respectively. The gradients of the two regression lines represent the strains, which are characterized by the gradients in the parenchyma, and show 1) higher gradients for seven out of ten subjects, 2) slightly higher gradients for Case 3, and 3)  lower gradients for the other two subjects.

Table \ref{table:strain} shows the mean strain for the bronchi and parenchyma region in the ten subjects. The mean values were computed from the strains of the sampled points in Fig. \ref{fig:strains}. The means $\pm$ standard deviation of the strains in the bronchus and lung parenchyma were 0.292 $\pm$ 0.141 and 0.395 $\pm$ 2.96, respectively. The strain of the lung parenchyma was higher than that of the bronchi, except for Cases 1 and 8. A significant difference was found between the two regions (ANOVA; $p < 0.05$ significance level), with the mean value of the parenchyma region being 35\% higher.

\begin{figure*}[t]
  \begin{center}
  \includegraphics[width=160mm]{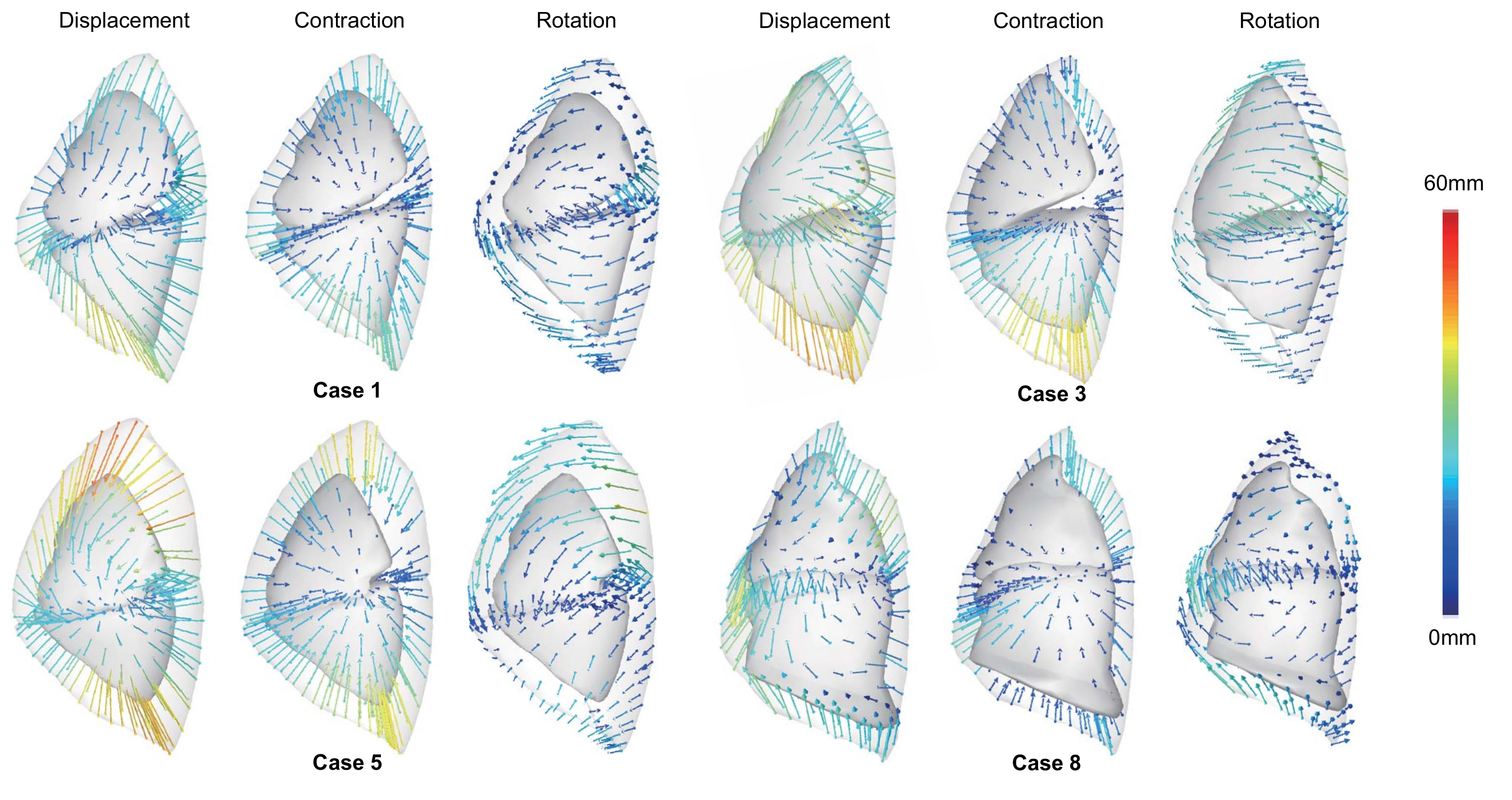}
  \caption{Visualization results of four deformation fields, contractions, and rotational components. The displacement vectors represent the corresponding points of the inflated/deflated models. The colors show the magnitude of the displacements. }
  \label{fig:rotations}
  \end{center}
\end{figure*}

We visualized the deformation fields to further investigate the reasons why such different characteristics appeared in the strains. Fig. \ref{fig:rotations} shows the visualization results of the displacements, contractions, and rotational components obtained from Cases 1, 3, 5, and 8. In Cases 1 and 8, the strain in the parenchyma is lower than that in the bronchi. In Case 3, the strains of the two regions are similar. Case 5 is typical, with a large deformation and strain in the parenchyma higher than that in bronchi. The displacement vectors represent the corresponding points of the inflated/deflated models and the colors show the magnitude of the displacements. 

The visualization results demonstrate the following points:

\begin{itemize}
    \item Rotational directions in upper and lower lobes were different, resulting in physical interactions and the occurrence of spaces among the lobes in Cases 1 and 3.
    \item In Case 5, both the upper and lower lobes deformed while preserving their shape features and relative postures.
    \item Deformation of the lower lobe in Case 8 was unique and was strongly affected by external forces from the diaphragm. 
\end{itemize}

We further discuss these influences on the strain analysis in the next section.

\section{Discussion}

% Contribution 
To our knowledge, this study is the first to analyze pneumothorax-associated deformation of whole lungs (surface and subsurface structures) using paired CT images measured in the same posture. Deformable registration methods were designed for inflated/deflated lung models, and non-linear mapping between the lung parenchyma and bronchi was obtained. While other studies have targeted pneumothorax-associated deformation, the main focus of their analyses was limited to either the lung surface or the bronchus region. A recent study \cite{Jeanne20} used CT data of six patients with pneumothorax after lung biopsy and chest tube insertion, and reported preliminary results for predicting the performance of deflated lung surfaces using a hyperelasitic finite-element model. Our previous study \cite{Nakao19} also focused on surface deformation of collapsed lungs.

% Strain analysis and irregular deformations
In the deformation analysis, strains in each region of the lung parenchyma and bronchi were calculated and a significant difference was found between the two strains. The mean strain of the parenchyma was higher than that of the bronchi, although considerable deviations in the values were found between subjects. In Case 3, the two regions showed similar strains, whereas Cases 1 and 8 showed the opposite trend, with the strain in the bronchi being higher than that in the parenchyma. The visualization results of the deformation fields in Fig. \ref{fig:rotations} demonstrated that Cases 1 and 3 showed opposite rotational directions between the lobes. In this situation, the two lobes made contact with each other and physical interaction could not be ignored. Although rotational components were found in other subjects, as shown in Case 5, the interaction between the lobes was observed to be small. Case 8 showed a unique deformation that was strongly affected by the diaphragm shape. We could not find the cause of this from the CT images, but soft tissue adhesion or biased forces from a ventilator might have acted on the lower lobe. The widely dispersed 2D plots in Fig. \ref{fig:rotations} suggest that the deformation in Case 8 contains outlying values with very small deformations or almost fixed points. We consider that these physical conditions influenced the natural contraction characteristics resulting in an error factor in the strain analysis of the three cases. 

% Limitation (only left lung, number of dataset)
This paper targeted the left lung to measure stable pneumothorax-associated deformation in a limited number of live dogs. In the right lung, the physical interaction between the three lobes (upper, middle, and lower) may be more complicated than that in the left lung, although the same measurement protocol and registration algorithms could be applied to the right lung. In our experiments, the imaging data were only collected from 11 subjects, and further acquisitions were difficult because of renovations being performed to our animal experiment facilities. Further analysis of the non-linear characteristics of both left and right deformed lungs would be interesting. 

% Clinical application and CBCT images
We assume that the clinical application of the proposed DMR techniques would be construction of a statistical deformation model representing the various pneumothorax states of the lung. Non-invasive intraoperative localization of tumors or nodules, especially those located deeper in the lungs, is desirable. The statistical deformation model could be used as prior knowledge to predict subsurface deformation in video-assisted thoracoscopic surgery \cite{Wu19}. Because CT imaging is not generally available during surgery, it will be challenging to construct a patient-specific image database of collapsed lungs in the intraoperative condition. However, we have begun to analyze intraoperative deformation using CBCT images \cite{Maekawa20}. Despite the limited measurement area, CBCT imaging is clinically feasible and will be useful for modeling intraoperative pneumothorax-associated deformations in real patients.

\section{Conclusion}
In this paper, we introduced a complete model-based registration solution for mapping the lung’s surface and bronchial structures including missing shapes, and analyzed the heterogeneity of pneumothorax-associated deformations using ten paired CT images of the lungs of Beagle dogs. The proposed DMR framework achieved stable registration with a Hausdorff distance error of less than 1 mm and a TRE error of less than 5 mm. Our results show that subsurface deformation of the collapsed lung was heterogeneous. The means $\pm$ standard deviation of the strains in the bronchi and lung parenchyma were 0.292 $\pm$ 0.141 and 0.395 $\pm$ 2.96, respectively, with the mean value of the parenchyma region being 35\% higher. We believe that the registered models and findings are suitable for statistical lung modeling and intraoperative tumor localization.

\section*{Acknowledgments} This research was supported by a JSPS Grant-in-Aid for Scientific Research (B) (Grant number 19H04484). Part of this study was also supported by a JSPS Grant-in-Aid for challenging Exploratory Research (Grant number 18K19918). We thank Karl Embleton, PhD, from Edanz Group (https://en-author-services.edanzgroup.com/ac) for editing a draft of this manuscript.

\bibliographystyle{IEEEtran}

\end{document}